\documentclass[10pt,reqno]{amsart}
\usepackage[centertags]{amsmath}
\usepackage{times}
\usepackage{mathptm}
\usepackage{mathrsfs}

\theoremstyle{plain} 

\newtheorem{lem}{Lemma}

\theoremstyle{definition}

\theoremstyle{remark}
\newtheorem{rem}{Remark}  

\numberwithin{equation}{section}

\DeclareMathSymbol{\R}{\mathalpha}{AMSb}{"52}
\DeclareMathSymbol{\C}{\mathalpha}{AMSb}{"43}

\newcommand{\comment}[1]{}
\newcommand{\set}[1]{\left\{#1\right\}}
\newcommand{\pd}{\,\partial}

\begin{document}

\title[]{Group Classification of a family of second-order differential equations}
\author[]{J.C. Ndogmo}

\address{PO Box 2446\\
Bellville 7535\\
South Africa }
\email{ndogmoj@yahoo.com}
\begin{abstract} We find the group of equivalence transformations
for equations of the form $y''= A(x)y' + F(y),$ where $A$ and $F$
are arbitrary functions.  We then give a complete group
classification of these families of equations, using a direct method
of analysis, together with the equivalence transformations.
\end{abstract}

%
\maketitle
\section{Introduction}

\label{s:intro} The group classification problem of equations of the
form
\begin{equation}\label{eq:2ode}
y''= F (x, y, y')
\end{equation}
was first considered by Lie \cite{lie-ovsi1}, who showed that the
symmetry group of all these equations is at most eight-dimensional
and that this maximum is reached only if the equation can be mapped
by a point transformation to a second-order linear ordinary
differential equation (ODE). More recently, Ovsyannikov considered
in \cite{ovsi1} the problem of group classification of a much
restricted form of the equation
 considered by Lie, namely the equation of the form
\begin{equation}\label{eq:ovsi1}
y''= F(x,y).
\end{equation}
 This study showed, amongst others, that in the
nonlinear case, which occurs if and only if $F=F(x,y)$ is not linear
in $y,$ the symmetry algebra has dimension at most two, except when
$F$  in \eqref{eq:ovsi1} can be reduced  by  an equivalence
transformation to $F= \pm y^{-3},$ in which case the nonlinear
equation has a symmetry algebra of maximal dimension three.\par 
Equations of the form \eqref{eq:2ode} containing a term linear in
$y',$ and of the form
\begin{equation} \label{eq:emden-f}
y''= \frac{M}{x} y' + g(x) F(y),
\end{equation}
where $g(x)$  and $F(y)$  are some given functions of $x$ and $y$
respectively, and $M$ is a constant, also appear frequently in the
mathematical physics literature. Eq. \eqref{eq:emden-f} is referred
to as an Emden-Fowler type equation \cite{waz1-c1, waz1-c2}, and it
is reduced for $M=2,\; g(x)=1$, and $F(y)= y^{\,n}$ to the so-called
{ \em standard} Lane-Emden equation of index $n,$ proposed by Lane
\cite{lane} and studied in detail by Emden \cite{emden} and Fowler
\cite{fowler}. It has been used as a model for the dynamics of a
spherical cloud of gas acting under mutual attraction of its
molecules \cite{waz1-c1}. Equation ~\eqref{eq:emden-f} with $g(x)=1$
is usually called the generalized Lane-Emden equation and for
special cases of $F(y),$ it has also been used as  a model for
various phenomena in physics and astrophysics, such as the  stellar
structure, the thermionic currents, and the dynamics of isothermal
spheres \cite{waz1, waz1-c1, waz1-c2, waz1-c3}.\par

Special cases as well as slightly modified forms of
\eqref{eq:emden-f} have been considered for symmetry analysis and
first integrals or exact solutions \cite{faz-c18, faz-c26, faz-c32}.
However, as far as the group classification of the Emden-Fowler type
equation is concerned, only a classification of Noether point
symmetries of the generalized Lane-Emden equation has been
considered \cite{fazlane}, but only for various (and not arbitrary )
functions $F.$  It should also be noted that the problem of
determination of solutions of \eqref{eq:emden-f} by analytic
approximations using the Adomian decomposition method, and
incorporating a singularity analysis was considered in ~\cite{waz1,
waz2}.\par
  The purpose of this paper is to provide a group classification of
the equation
\begin{equation}\label{eq:m}
y''= A(x)y'+F(y),
\end{equation}
in which $A$ and $F$ are arbitrary functions of the independent
variable $x$ and the dependent variable $y,$ respectively. This is a
modified form of \eqref{eq:emden-f} in which the coefficient of $y'$
is an arbitrary function and $g(x)=1.$ We find the group of
equivalence transformations of this equation, that is, the largest
group of point transformations that preserves the form of the
equation. Next, we obtain the group classification of the equation,
based on a direct analysis and using the equivalence
transformations. It is shown in particular that in the nonlinear
case, which occurs if and only if $F$ is nonlinear, the maximal
dimension of the symmetry algebra is three. Moreover, it is also
clearly shown that any symmetry exists only for canonical forms of
$F$ of the form
\begin{equation}\label{eq:admis}
f,\quad  \mu e^y + f, \quad \mu y \ln (y) +f,\quad \; \mu \ln(y)
+f,\quad \text{ and } y^n +f,\quad  n\neq 0,1,
\end{equation}
where $\mu\neq0$ is a constant and $f$ is a linear function of $y.$

\section{Equivalence group}

We shall say that an invertible point transformation of the form
\begin{equation}\label{eq:eqvtf}
x= S(z,w), \qquad y= T(z,w)
\end{equation}
is an equivalence transformation of \eqref{eq:m} if it transforms
the latter equation into an equation of the same form, that is, into
an equation of the form
\begin{equation}\label{eq:mtf}
w''= B(z) w\,' + H(w),
\end{equation}
where $B(z)$ and $H(w)$ are the new arbitrary functions, and where
$w\,'= ~d w/ dz.$ In this case the two equations \eqref{eq:m} and
\eqref{eq:mtf} are said to be equivalent. The equivalence group $G$
of \eqref{eq:m} is the largest Lie pseudo-group of transformations
of the form \eqref{eq:eqvtf} that preserves the form of the
equation. By a result of Lie \cite{lieGc}, the resulting
transformations of the arbitrary functions $A$ and $F$ also form a
Lie pseudo-group of transformations, which in the actual case can be
put into the form
\begin{equation}\label{eq:coftf}
A= \chi(z, w, B, H), \quad \qquad F= \zeta (z, w, B, H),
\end{equation}
for certain functions  $\chi$ and $\zeta$ which may be read-off from
 expressions of the transformed coefficients once the defining
functions $S$ and $T$ of $G$ are known. By writing down the
transformation of \eqref{eq:m} under \eqref{eq:eqvtf}, we obtain an
equation in  $S$ and $T$  which is rearranged by an expansion into
powers of $w'$. Then, using the fact that $S$ and $T$ do not depend
explicitly on the derivatives of $w$ with respect to $z,$ and
assuming that \eqref{eq:eqvtf}  maps \eqref{eq:m} to \eqref{eq:mtf},
the transformed equation is reduced to the following set of four
equations in which $B(z) w' + H(w)$ is substituted for $w'',$ and
where $\delta= \left(S_z T_w-S_w T_z\right)$ is the Jacobian of the
change of variables \eqref{eq:eqvtf} $\colon$
\begin{subequations}\label{eq:deteqv}
\begin{align}
  &-F S_z^3+H \delta -T_z \left(A S_z^2+S_{z,z}\right)+S_z T_{z,z}=0 \label{eq:deteqv1}\\
 \begin{split}  &-3 F S_w S_z^2+B S_z T_w-A S_z^2 T_w-B S_w T_z-2 A S_w S_z T_z\\
&-2 T_z S_{z,w}-T_w S_{z,z}+2 S_z T_{z,w}+S_w T_{z,z}=0 \end{split}\label{eq:deteqv2}\\
 &-3 F S_w^2 S_z-A S_w \left(2 S_z T_w+S_w T_z\right)-T_z S_{w,w}-2 T_w S_{z,w}+S_z T_{w,w}+2 S_w T_{z,w}=0 \label{eq:deteqv3}\\
 &-F S_w^3-T_w \left(A S_w^2+S_{w,w}\right)+S_w T_{w,w} =0. \label{eq:deteqv4}
\end{align}
\end{subequations}
From \eqref{eq:deteqv4}, it follows that $S_w^3=0,$ on account of
the arbitrariness of $F,$ and so $S=S(z).$ When this last equality
is substituted into \eqref{eq:deteqv}, \eqref{eq:deteqv3} is reduced
to $S_z T_{w,w}=0,$ which shows that $T= \alpha(z)w+ \beta(z)$, for
some functions $\alpha$ and $\beta.$ With these expressions for $S$
and $T,$ the first two equations of \eqref{eq:deteqv} are reduced to
\begin{subequations}\label{eq:rdeqv}
\begin{align}
 &H \alpha  S_z-F S_z^3-A S_z^2 \left(w \alpha _z+\beta _z\right)-
\left(w \alpha _z+\beta _z\right) S_{z,z}+S_z \left(w \alpha
_{z,z}+\beta
 _{z,z}\right)=0 \label{eq:rdeqv1}\\
 &B \alpha  S_z-A \alpha  S_z^2+2 S_z \alpha _z-\alpha  S_{z,z}=0. \label{eq:rdeqv2}
\end{align}
\end{subequations}

In \eqref{eq:rdeqv1},  $H$ depends {\em a priori} on both $F$ and
$A,$ but in virtue of the conditions $H=H(w)$ and $A=A(z),$
 $H$,  (and clearly F), must be independent of $A.$ Consequently, it follows
from the arbitrariness of $A$ that its coefficient $(w \alpha
_z+\beta _z)$ in \eqref{eq:rdeqv1} must vanish identically.
Therefore, $\alpha= k_3$ and $\beta= k_4,$ for some constants  $k_3$
and $k_4.$ Substituting these values for $\alpha$ and $\beta$ into
\eqref{eq:rdeqv1} and solving for $H$ gives
\begin{equation}\label{eq:eq4h1}
H= F(k_3+ k_4 w) S_z^2/ k_3,
\end{equation}
and the condition $H_z=0$ forces $S_z$ to be a constant function.
Consequently, we must have $S= k_1 z+ k_2,$ where $k_1$ and $k_2$
are arbitrary constants. A substitution of the expressions thus
obtained for $\alpha, \beta$  and $S$ into \eqref{eq:rdeqv}
completely determines $B$ and $H.$   The group $G$ of equivalence
transformations of \eqref{eq:m}  is therefore given by the linear
transformations
\begin{equation}\label{eq:equiv1}
x= k_1 z+ k_2, \qquad y = k_3 w + k_4, \quad (k_1 k_3= \delta \neq
0),
\end{equation}
where the $k_j,$ for $j=1, \dots, 4$ are constants. On the other
hand, the resulting induced transformations of the arbitrary
functions $A$ and $F$ are given by
\begin{equation}\label{eq:equiv2}
B= k_1 \,A \circ L^{k_1, k_2}, \qquad H= \frac{k_1^2}{k_3} \,F \circ
L^{k_3, k_4}
\end{equation}
where  for every pair of constants $p$ and $q,$ $L^{\,p,q}$ is the
linear function $L^{\,p,q} (a)= p\,a+q.$ Consequently, the explicit
form of the transformed equation \eqref{eq:mtf} under
\eqref{eq:equiv1} is given by
\begin{equation}\label{eq:m2tf}
w''= k_1 A(k_1 z + k_2) + \frac{k_1^2}{k_3} F(k_3 w + k_4),
\end{equation}
and this shows in particular that for arbitrary values of the
functions $A$ and $F,$ \eqref{eq:m} has no nontrivial symmetries.
Indeed, the symmetry group of the equation is a subgroup of its
equivalence group, and \eqref{eq:m2tf} shows that \eqref{eq:equiv1}
is a symmetry of \eqref{eq:m} for every functions $A$ and $F$ only
if
\begin{equation*}k_1 z + k_2= z, \qquad \text{ and} \qquad k_3 w + k_4=w,
\end{equation*} that is, only if \eqref{eq:equiv1} is the identity
transformation.\par
An element of the family of equations of the form \eqref{eq:m} may
be labeled by the corresponding pair $\set{A, F}$ of coefficient
functions, and by the above  results two equations represented by
the pairs $\set{A, F}$ and $\set{B,H}$ are equivalent under
\eqref{eq:eqvtf} if the coefficient functions are related by
\eqref{eq:equiv2}. By changing only the dependent variable in
\eqref{eq:m}, we may keep $A$ fixed and transform only  $F,$ which
induces among the coefficient functions $F$ another equivalence
relation that we denote by $\sim.$  We have the following result
about the latter equivalence relation.
\begin{lem}
Let $a,b, c, r, s$ and  $n$ be given constants with $a, r\neq0$ and
$n\neq 1.$ There are constants $\mu \neq0,$ and $\lambda, \theta$
such that the following hold in each case.
\begin{enumerate}
\item[{\rm (a)}] $r(ay +b)^n + c y + s$ $\sim$ $y^n + \lambda y + \theta$,
and $a y^2 + b y + c$ $\sim$ $y^2 + \theta.$
\item[{\rm (b)}] Let $F= r\, e^{a\, y} + b y + c.$ Then $F \sim \mu e^y + \lambda
y$ for $\lambda \neq 0.$ Else $F \sim \mu\, e^y + \theta.$
\item[{\rm (c)}] $a \ln (y) + b y + c$ $\sim$ $\mu\, \ln(y) + \lambda y.$
\item[{\rm (d)}] $a y \ln (y) + b y + c$ $\sim$ $\mu\, y \ln(y) + \theta.$
\item[{\rm (e)}] Let $F= c y + b.$ Then $F \sim \mu y$ for $c\neq0,$ else
$F\sim \theta,$ where $\theta=0$ or $\theta=1.$
\end{enumerate}
\end{lem}
\begin{proof}
According to \eqref{eq:equiv2}, we only need to show that in each
case we can find  constants $k_3$ and $k_4$ such that the given
function $F$ is equivalent to the indicated function $H$ for some
constants $\mu,$ $\lambda$ and $\theta$ to be specified. This is
achieved by finding for the given function, $F(y)$ say, the
transformed function  $H= F(k_3 y + k_4)/k_3$ which has the required
form. In case (a) for example, letting
\begin{equation*}F=r(ay+ b)^n + c y + s,\end{equation*}
it is readily seen that by choosing $k_3= (r\, a^n)^{1/(1-n)}$ and
$k_4= -b/a,$ $F$ is transformed into $H=y^{\,n} + \lambda y +
\theta,$ with $\lambda =c$ and $\theta= -(b c) / (a k_3) + s/k_3.$
For the second part of (a) with $F= a y^2 + b y + c,$ the result
follows by setting $k_3= 1/a,$ $k_4= -b/(2a),$ which gives $\theta=
(4 ac - b^2)/4.$ The other cases are treated in a similar manner.
\end{proof}
\section{Group classification}
Equivalence transformations are often very helpful for the group
classification of differential equations, because equivalent
equations also have equivalent symmetry algebras, in the sense that
one can be mapped onto the other by an invertible change of
variables. However, in the actual case of Eq. \eqref{eq:m}, the
transformations obtained in \eqref{eq:equiv1} and \eqref{eq:equiv2}
are relatively weak, in the sense that they act on arbitrary
functions only by mere scalings and translations and give rise in
particular to an infinity of non equivalent equations. Nevertheless,
we shall still be able to use them as a simplifying tool in our
classification procedure of \eqref{eq:m}, which is based on a direct
analysis of the determining equation of the symmetry algebra. In the
sequel, the symbols $M$ and $m$, as well as $k_1, k_2, \dots$ will
represent arbitrary constants. For a given function $Q=Q(a)$ with
argument $a,$ we shall write $Q\,'$ for $d Q/da.$  If we let
\begin{equation}\label{eq:Lgen}
    V= \xi (x,y) \pd_x + \phi(x,y)\pd_y
\end{equation}
denote the generic generator of the symmetry algebra $L$ of
\eqref{eq:m}, then it readily follows from well-known procedures
\cite{olvbk86, bluman} that the determining equations of $L$ are
given  for arbitrary functions $A$ and $F$ by
\begin{subequations}\label{eq:detm}
\begin{align}
 & \xi _{y,y}=0   \label{eq:detm1}\\
 &-\xi  A'-A \xi _x-3 F \xi _y-\xi _{x,x}+2 \phi _{x,y}=0  \label{eq:detm2}\\
 &-\phi  F'-2 F \xi _x-A \phi _x+F \phi _y+\phi _{x,x}=0 \label{eq:detm3}\\
 &-2 A \xi _y-2 \xi _{x,y}+\phi _{y,y}=0. \label{eq:detm4}
\end{align}
\end{subequations}
From \eqref{eq:detm1} and \eqref{eq:detm4}, it follows successively
that
\begin{equation}\label{eq:xiphi}
\xi= \alpha (x) y+ \beta(x), \qquad \text{ and } \qquad \phi= y^2(A
\alpha(x) + \alpha'(x)) + y\, \sigma(x) + \tau(x),
\end{equation}
for some functions $\alpha, \beta, \sigma$ and $\tau.$ Next, a
substitution of \eqref{eq:xiphi} into \eqref{eq:detm}  transforms
\eqref{eq:detm2} and \eqref{eq:detm3} into
\begin{subequations}\label{eq:red1detm}
\begin{align}
 &-3 F \alpha   +3 y \left(\alpha  A'+A \alpha'+\alpha''\right)
 -\beta  A'-A \beta\,'+2 \sigma \,'-\beta\,''=0 \label{eq:red1detm1}\\
 \begin{split}&-F' \left(-y\, \sigma -\tau +y^2 \left(-A \alpha -\alpha'\right)\right)
 +F \left(2 A y\, \alpha +\sigma -2 \beta\,'\right)-A \tau\,' +\tau\,'' \\
 &+y \left(-A \sigma
\,'+\sigma\,''\right)+y^2 \left(-A \alpha  A'-A^2 \alpha'+2 A'
\alpha'+\alpha  A''+\alpha\,'''\right)=0 \label{req:ed1detm2}
\end{split}
\end{align}
\end{subequations}

Differentiating twice \eqref{eq:red1detm1} with respect to $y$ shows
that $-3 F''\, \alpha=0,$ and thus we shall consider separately the
cases $F''=0$ and $F'' \neq 0.$\par

\subsection{Case 1: $F'' \neq 0$}
We must have in this case $\alpha=0,$ and  this reduces
\eqref{eq:red1detm} to
\begin{subequations}\label{eq:alpha0}
\begin{align}
& -\beta  A'-A \beta\,'+2 \sigma \,'-\beta\,''=0 \label{eq:alpha01}\\
& (-y \sigma -\tau ) F'+F \left(\sigma -2 \beta\,'\right)+y \left(-A
\sigma \,'+\sigma\,''\right)  -A \tau\,'+\tau\,'' =0
\label{eq:alpha02}
\end{align}
\end{subequations}
Differentiating \eqref{eq:alpha02} with respect to $y$ twice yields
\begin{equation}\label{eq:valF1}
    (\sigma + 2 \beta\,')F'' + (y\, \sigma + \tau) F'''=0.
\end{equation}
If $F'''=0,$ then by the lemma we may assume that $F = y^2 + \theta$
where $\theta$ is a constant, and a substitution of the latter
expression for $F$ into \eqref{eq:alpha02} leads after an expansion
into powers of $y$ to the equality
\begin{equation}\label{eq:sigma}
\sigma = -2 \beta\,',
\end{equation}
representing the vanishing of the coefficient of $y^2.$ Substituting
also the latter expression for $\sigma $ into \eqref{eq:alpha0} and
expanding again into powers of $y$ gives
\begin{subequations}\label{eq:alpha022}
\begin{align}
-\beta  A'-A \beta\,'-5 \beta\,''&=0 \label{eq:alpha022a}\\
-4 \theta\, \beta\,'-A \tau \,'+\tau\,'' &=0 \label{eq:split1a}\\
 \tau -A \beta\,''+\beta\,'''&=0, \label{eq:split1b}
\end{align}
\end{subequations}
\begin{rem}\mbox{$\quad$}
\begin{enumerate}
\item[(a)] It appears from the total orders of derivatives of unknown  functions
appearing in \eqref{eq:alpha022} that the maximal possible number
$\kappa_M$ of free parameters in the general solution cannot exceed
five and therefore five is an upper bound for the dimension of the
corresponding symmetry algebra $L.$ However, the constraints in the
system will in general reduce the size of $\kappa_M.$
\item[(b)] The order in which individual equations are solved and
the corresponding solutions are substituted into the system doesn't
matter, for they all yield the same solution. For
\eqref{eq:alpha022} and all subsequent similar systems of equations,
we shall therefore choose the order of integration that appears to
be the most suitable for the solution.
\end{enumerate}
\end{rem}

For \eqref{eq:alpha022}, a suitable order of integration consists in
 solving \eqref{eq:alpha022a} for $\beta,$ substituting the result
into \eqref{eq:split1b} to find $\tau$, and using \eqref{eq:split1a}
for the resulting compatibility condition on $A.$ This, together
with \eqref{eq:sigma}, yields
\begin{subequations}\label{eq:gen1}
\begin{align}
\beta &=  e^{ -\int \frac{A}{5}
 \, dx} \left(k_2+ k_1 \int e^{
\int \frac{A}{5} \, dx} \, dx \right) \label{eq:gen1a}\\
\sigma &= -2 k_1+\frac{2}{5} e^{-\int \frac{A}{5} \, dx} A \left(k_2+k_1 \int e^{ \int \frac{A}{5} \, dx} \, dx\right) \label{eq:gen1b}\\
\tau   &= \frac{1}{125} \left[-30 k_1 A^2+50 k_1 A'+e^{- \int
\frac{A}{5} \, dx} \left(k_2+k_1 \int e^{ \int
 \frac{A}{5}
\, dx} \, dx\right) \left(6 A^3-40 A A'+25 A''\right)\right],
\label{eq:gen1c}
\end{align}
\end{subequations}

and the compatibility  condition on the coefficient $A$ for the
existence of any symmetry is given by
\begin{subequations}\label{eq:condac1}
\begin{align}
k_1 &(- 20 F_2 E_2+ F_1 E_1 ) + k_2 E_1=0 \label{eq:condac1a}\\
\intertext{ where }
 \begin{split} E_1 &= 36 A^5-900 A^3 A'+2000 A^2 A''\\
        &\quad +625 A \left(4 \left(A'^2+\theta\right) -3 A^{(3)}\right)+625 \left(-5 A'\, A''+A^{(4)}\right)\end{split} \label{eq:condac1b}\\
E_2 &= 9 A^4-180 A^2 A'+275 A A''+25 \left(7 A'^2+25 \theta-5 A^{(3)}\right)  \label{eq:condac1c}\\
F_1 &=\int e^{\frac{1}{5} \int A \, dx} \, dx, \qquad \text{ and
}\quad  F_2= e^{ -\int A/5 \, dx}. \label{eq:condac1d}
\end{align}
\end{subequations}
It is worthwhile recalling that the expression of the symmetry
generator $V$ of \eqref{eq:Lgen} is reduced in this case to
\begin{equation*}
V= \beta(x)\pd_x + \left( y\, \sigma(x) + \tau(x) \right) \pd_y.
\end{equation*}
 For an explicit determination of the dimension of $L,$ we note that since
$E_1= -5 E_2'+ 4 A E_2,$ where $E_2'= d E_2/dx,$ it follows that
$E_1=0$ if $E_2=0,$ and hence $L$ has dimension two if and only if
$E_2=0.$  On the other hand, $L$ has dimension one if and only if
exactly one of the following conditions hold
\begin{equation}\label{eq:c0y2dim1}
  E_1=0, \qquad \text{ or } - 20 F_2 E_2+ F_1 E_1=0,
\end{equation}
the latter condition being an integro-differential equation.
One-parameters families of solutions of \eqref{eq:condac1a} indexed
the arbitrary constant $m$ are given by
\begin{subequations}\label{eq:sola4y2}
\begin{align}
A &= p / (x+m), \quad &\text{for } \theta &=0, \text{ and }p=0,-15, -10/3, -5/3 \label{eq:sola4y2a}\\
A &= 5 p \tan(p x+m), \quad & \text{for } p &=  (- \sqrt{\theta}\, i
/3)^{1/2}, \text{ and $k_1=0,$} \label{eq:sola4y2b}
\end{align}
\end{subequations}
and in case \eqref{eq:sola4y2a} $L$ has dimension two, while it has
dimension one in case \eqref{eq:sola4y2b}, with generator
\begin{equation}\label{eq:geny2nz}
\cos \left[m+\frac{x \sqrt{-i \sqrt{\theta
 }}}{\sqrt{3}}\right]\pd_x +
 \frac{2 \left(y+i \sqrt{\theta }\right) \sqrt{-i \sqrt{\theta }} \sin \left[m+\frac{x \sqrt{-i \sqrt{\theta
 }}}{\sqrt{3}}\right]}{\sqrt{3}}\pd_y.
\end{equation}
We also have $\dim L=1$ for $A=M,$ but we can  hardly describe all
possibilities when $L$ has exactly dimension one or  two according
to the values of $A,$ because general solutions of $E_2=0,$ where
$E_2$ is given by \eqref{eq:condac1c}, or of \eqref{eq:c0y2dim1}
aren't available. This completes the classification problem when
$F'''=0.$\par
   If $F''' \neq 0,$  it follows from
\eqref{eq:valF1} that $(F''/F''')''=0,$ and hence
\begin{equation*}F'''/F'' = 1/ (a_1 y + a_2)
\end{equation*}
for some constants $a_1$ and $a_2,$ and  we have to consider
separately  the two possibilities  $a_1 \neq 0,$ and $a_1=0$ (but
$a_2 \neq 0).$  All these possibilities lead to the following
possible canonical forms for $F:$
\begin{subequations}\label{eq:F-forms}
\begin{alignat}{4}
&\text{Case (i): } &F &= \mu \, e^y + \lambda y \quad (\lambda
\,\neq 0), \qquad &\text{Case (ii):}\quad & F&=
\mu \, e^y + \theta\\
&\text{Case (iii): } &F &= \mu \, \ln (y) + \lambda y,  \qquad
&\text{Case (iv):}\quad & F &= \mu \, y \ln (y) +
\theta\\
&\text{Case (v): } & F& = y^{\,n} + \lambda y + \theta,\; (n \neq
0,1,2),\quad & && &
\end{alignat}
\end{subequations}
where $\lambda, \theta$ and $\mu$ are constants with $\mu \neq 0$,
and where the first two cases (i) and (ii) correspond to $a_1=0,$
while the remaining cases correspond to $a_1 \neq 0$.\par
   In the cases (i), (iii), and the cases (iv) and (v) with $\theta
\neq 0,$ it is readily found that the only symmetry is $V= \pd_x,$
provided that the coefficient $A$ is a constant function.\par

In case (ii), when $\theta=0$ we have $F= \mu e^y$  and a
substitution  of this expression into \eqref{eq:alpha0} shows that
$L$ is generated by
\begin{equation}\label{eq:theta0}
V= \left(k_2 + k_3 x + k_1 \int \int e^{\int A dx} dx dx\right)
\pd_x + \left(-2 k_3 - 2 k_1 \int e^{\int A dx}dx\right)\pd_y
\end{equation}
provided that the compatibility condition
\begin{equation}\label{eq:condac2}
-k_2 A'- k_3 \left(A+x A'\right)+k_1 \left(-e^{\int A\, dx}-A\int
e^{\int A\, dx} \, dx-  A' \int \int e^{\int A\, dx} \, dx \, dx
\right)=0
\end{equation}

on the coefficient $A$ is satisfied. An analysis of
\eqref{eq:condac2} shows that $\dim L =2$ for $A=0$ or $A=-1/x,$ and
$\dim L=1$ for $A=M$, $M\neq 0$ or $A=M/x,$ $M \neq -1.$ Else $L$ is
zero-dimensional.\par
When $F= \mu \, e^y + \theta,$ and $\theta \neq 0,$ a substitution
of this expression into \eqref{eq:alpha0} shows that we must have
$\alpha=\sigma=0,$ and $\tau= -2 \beta'.$ All these expressions for
$F, \alpha, \sigma$ and $\tau$ reduce \eqref{eq:alpha0} to
\begin{subequations}\label{eq:theta1}
\begin{align}
- \beta A' - A \beta' -\beta''& =0  \label{eq:theta1a}\\
\theta\, \beta' - A \beta'' + \beta'''&=0 \label{eq:theta1b},
\end{align}
\end{subequations}
and we find the solution by  solving \eqref{eq:theta1a}
 for $\beta$ and substituting the result into the second equation. In this way we find the
corresponding symmetry generator to be of the form
\begin{subequations}\label{eq:gentheta}
\begin{align}
V &= F_2  (k_2 + k_1 F_1)\pd_x + \left(-2 k_1 +
2F_2 A (k_2 + k_1 F_1)   \right) \pd_y \label{eq:gentheta1}\\
\intertext{ where }
F_1 &= \int e^{\int A \, dx} \, dx, \qquad F_2= e^{-\int A \, dx}
\label{eq:gentheta2}
\end{align}
\end{subequations}
and the compatibility condition on $A$ is given by
\begin{subequations}\label{eq:condac3}
\begin{align}
 k_1 & \left(- E_4 + F_2\, F_1 E_3\right) +  k_2
F_2 E_3 =0  \\
\intertext{where}
E_3 &=2 A^3+A \left(\theta-4 A'\right)+A'',\quad  \text{ and }\quad
E_4= \theta+2 A^2-2 A'.
\end{align}
\end{subequations}
We have $E_4' + 2 E_3= 2 A E_4$ where $E_4' = d E_4/dx,$ and thus we
readily see that $L$ has dimension two if and only if $E_4=0,$ that
is, if and only if
\begin{equation*} A= \sqrt{\frac{\theta}{2}} \tan \left( \sqrt{\frac{\theta}{2}}  (x+ 2m)
\right),
\end{equation*}
where $m$ is a constant. On the other hand,  $L$ has dimension one
if $A$ satisfies exactly one of the conditions
\begin{equation}\label{eq:c2nzdim1}
 E_3=0, \quad \text{ or } - E_4 + F_1 F_2 E_3=0.
\end{equation}
 More explicitly, the latter
equality is an integro-differential equation of the form
\begin{equation}\label{eq:c2nzdim1b}
    -2 \left(\theta +2 A^2-2 A_x\right)+2 e^{-\int A \, dx} \left(\int e^{\int A \, dx} \, dx\right) \left(2 A^3+A \left(\theta -4
    A_x\right)+A_{x,x}\right)=0.
\end{equation}
In case {\rm (iv)} with $\theta = 0,$ we have $F= \mu\, y \ln(y),$
and a substitution of this expression for $F$ into \eqref{eq:alpha0}
shows that  $\alpha= \tau=0,$ and $\beta=k_1,$ and the remaining
conditions on $A$ and $\sigma$ are given by
\begin{subequations}\label{eq:caseiv}
\begin{align}
- k_1 A' + 2 \sigma'& =0  \label{eq:caseiva}\\
- \mu\, \sigma\, - A \sigma' + \sigma''&=0, \label{eq:caseivb}
\end{align}
\end{subequations}
Since for every function $A$ \eqref{eq:caseiv} always consists of a
first-order and a second-order ODE, its solution will depend on at
most one arbitrary constant, and $L$ will therefore have at most
dimension two in this case, with corresponding symmetry generator
\begin{equation*} V= k_1 \pd_x + y\, \sigma \, \pd_y.\end{equation*}
If we solve \eqref{eq:caseiva} for either $A$ or $\sigma$ and
substitute the result in \eqref{eq:caseivb}, then neither the
resulting equation, nor \eqref{eq:caseivb} itself is tractable. It
is therefore not possible to describe all the possible dimensions of
$L$ according to the values of $A.$ It is however clear that for
$A=M,$ we have $\sigma=0,$ and $V= \pd_x.$\  More generally, the
compatibility condition on $A$ is given by
\begin{equation}\label{eq:condac4z}
2 k_2 \mu + k_1 A (\mu +A')- k_1 A''=0
\end{equation}
where $k_2$ is another constant. \par
For case {\rm (v)}, when $\theta=\lambda=0,$ we have $F= y^ {\,n}$
$(n\neq 0,1,2),$ and the corresponding generator of $L$ is given by
\begin{equation*}
V= \left(k_2 + k_3 x + k_1 \int \int e^{\int A dx} dx dx\right)
\pd_x - \left(\frac{2 k_3 y} {n-1} + \frac{2 k_1 y \int e^{\int A
dx} dx }{n-1} \right) \pd_y
\end{equation*}
while the related compatibility condition on $A$ is given by
\begin{equation}
\begin{split}
&-k_2 (n-1) A'-k_3 (n-1) \left(A+x A'\right)\\
&-k_1\left(e^{\int A \, dx} (3+n)+(n-1) A \int e^{\int A \, dx} \,
dx+(n-1) A'\int \int e^{\int A \, dx} \, dx \, dx\right)=0.
\end{split}
\end{equation}
It thus follows that $L$ has dimension three if and only if $A=0$
and $n=-3.$ On the other hand, we have $\dim L=2$ if
\begin{alignat*}{2}
A &= -\frac{(n+3)/x}{n+1 }, \quad n\neq -1, -3,&\quad  &\text{ or if }\\
A &=M, \quad \text{with}\quad  M\neq 0 \quad\text{and}\quad n=-1,&\quad  &\text{ or if }\quad\\
A & =0 \quad\text{and}\quad n\neq -3.& &
\end{alignat*}

 We also have $\dim L=1$ if
\begin{align*}
A=\begin{cases} M/x, & \text{ with } M  \neq 0, -(3+n)/(n+1)\\
 M, & \text{ with }  M  \neq 0,\quad  \text{and}\quad
n\neq -1
\end{cases}
\end{align*}
or if
\begin{subequations}\label{eq:c5yndim1}
\begin{align}
A & \neq -\frac{(n+3)/x}{n+1 }, \quad   \text{ and }\\
0&= e^{\int A \, dx} (3+n)+(n-1) A \int e^{\int A \, dx} \, dx+(n-1)
A'\int \int e^{\int A \, dx} \, dx \, dx
\end{align}
\end{subequations}
 In case (v)  with  $\lambda \neq 0$, the symmetry generator is given by
\begin{equation}\label{eq:casevgen}
V= \beta \pd_x - \frac{2 y\, \beta'}{n-1}\pd_y,
\end{equation}

where $\beta$ is determined together with the coefficient $A$ by the
equation
\begin{subequations}\label{eq:detc5}
\begin{align}
\beta A' + A \beta' - (3+n) \beta'' &=0 \label{eq:detc5a}     \\
(n-1)\lambda  \beta' -A  \beta'' +  \beta\,''' &=0.
 \label{eq:detc5b}
\end{align}
\end{subequations}
 It follows from
\eqref{eq:detc5a} that for $n\neq -3$ we have

\begin{subequations}\label{eq:betac5}
\begin{align}
\beta &= (k_2 + k_1 F_1) F_2 \\
\intertext{where}
F_1&= \int \exp\left(\frac{(n-1) \int A \, dx}{3+n}\right) \, dx
\quad \text{
 and }\quad F_2= - \exp\left(\frac{(n-1) \int A \, dx}{3+n}\right)
\end{align}
\end{subequations}
and the corresponding compatibility conditions for $A$ are given by
\begin{subequations}\label{eq:condac5}
\begin{align}
k_2 &F_2 E_5 + k_1 \left((n+3) E_6 + F_1 F_2 E_5\right) =0,
\intertext{ where }
E_5  &= 2 A^3 \left(-1+n^2\right)+A (3+n) \left((-1+n) (3+n) \lambda -4 n A'\right)+(3+n)^2 A'' \\
E_6  &= -2 A^2 (1+n)+(3+n) \left((-3-n) \lambda +2 A'\right).
\end{align}
\end{subequations}
In this case we have
\begin{equation*}
2 E_5 - (3+n)E_6' + 2(n-1) A E_6=0,
\end{equation*}
where $E_6'= d E_6/dx.$ Consequently, we have $\dim L=2$ if and only
if
\begin{equation*}
A = \begin{cases}\frac{(3+n) \sqrt{\lambda } } { \sqrt{2(1+n)}} \tan
\left[\frac{\sqrt{\lambda(n+1)}\left(x+2
(3+n)m\right)}{\sqrt{2}}\right],&\quad
\text{for } n \neq -1\\
\lambda x + m,&\quad \text{for }n =-1
\end{cases}
\end{equation*}
where $m$ is a constant parameter. Clearly, we have $\dim L= 1$ if
exactly one of the following conditions hold:
\begin{equation}\label{eq:c5nzdim1}
E_5=0, \quad \text{ or } \quad (n+3) E_6 + F_1 F_2 E_5=0.
\end{equation}
However, we cannot obtain in general all explicit expressions of $A$
for which $\dim L=1,$ although here again for $A=M$ we have $\dim
L=1,$ and $V= \pd_x.$
\par
For $n=-3,$ Eq. \eqref{eq:detc5} reduces to
\begin{subequations}\label{eq:detc5m3}
\begin{align}
\beta A' + A \beta'  &=0 \label{eq:detc5m3a}     \\
-4 \lambda  \beta' -A  \beta'' +  \beta\,''' &=0.
 \label{eq:detc5m3b}
\end{align}
\end{subequations}
For $A=0,$ \eqref{eq:detc5m3a} vanishes identically, and
\eqref{eq:detc5m3} reduces to  \eqref{eq:detc5m3b} in which $A=0.$
~Solving the resulting equation for $\beta$ and substituting the
result into $\eqref{eq:casevgen}$ yields the generator
\begin{equation}\label{eq:vc5m3}
 V= \left[k_3 + \frac{e^{-2x \sqrt{\lambda}}\left( k_1 e^{4 x
\sqrt{\lambda}}  -k_2   \right) }{2\sqrt{\lambda}} \right]\pd_x +
\frac{1}{2}e^{-2x \sqrt{\lambda}}  \left( k_1 e^{4 x \sqrt{\lambda}}
+k_2   \right)y\pd_y
\end{equation}
where $k_1, k_2$ and $k_3$ are arbitrary constants and this shows in
particular that $L$ has dimension three in this case. For $A\neq 0,$
the general solution of \eqref{eq:detc5m3} will depend on at most
one arbitrary parameter, and hence $L$ will have dimension at most
one. If we set for instance  $A=M$ in \eqref{eq:detc5m3}, where $M$
is a nonzero constant, this gives $\beta= k_1,$ where $k_1$ is
another constant, with corresponding symmetry generator $V= k_1
\pd_x.$ However, we cannot describe all other  solutions to
\eqref{eq:detc5m3}, and hence we cannot describe in this case all
values of $A$ for which the dimension of $L$ takes on the value one.
Indeed, from \eqref{eq:detc5m3a} we have $A= k_1/\beta,$ and
substituting this into \eqref{eq:detc5m3b} gives
\begin{equation*}
-4 \lambda \beta' - k_1 \beta''/\beta + \beta^{(3)}=0,
\end{equation*}
which  is an equation for which the general solution is not
available. On the other hand we can look for one-dimensional
subalgebras of $L$ by solving  \eqref{eq:detc5m3a}
 for $\beta.$ This gives $\beta= k_1/A,$ and the corresponding
compatibility condition for $A$ takes the form
\begin{equation}\label{eq:condam3}
2 A'^2 + \frac{6 A'^3}{A^2} - A A'' + A'\left( - 4 \lambda - \frac{6
A''}{A} \right) + A'''=0.
\end{equation}
However, for this nonlinear equation we can only obtain the
particular solution $A= {const.}$ 
\begin{table}[h]

 \caption{\label{tb:summary} \protect \footnotesize  Classification results for the equation $y''= A(x)
y' + F(y).$} { \footnotesize
\begin{tabular}{l c c l}\hline\hline
 $F$  &   $A$      &     $\dim L$   & Generator $V$
 \\[1pt]\hline
    $\mu\, e^y $      &  $0$  &   2    & $(k_1 + k_2 x)\pd_x - 2k_2\pd_y$       \\[1pt]
                      &  $-1/x$ &   2    &       \\[1pt]
                    &   $M/x,\quad M\neq -1$      &   1    & $x \pd_x -2 \pd_y$
                    \\[1pt]\hline
$\mu\, e^y+ \theta, \; \theta \neq 0$ &$ \sqrt{\frac{\theta}{2}}
\tan \left( \sqrt{\frac{\theta}{2}}  (x+ 2m) \right) $ & 2&
\\[1pt]
          &    As given by \eqref{eq:c2nzdim1}    &  1
          &\\[1pt]\hline

 $\mu\, y \ln (y)$  &As given by \eqref{eq:condac4z}      &1 or 2& \\\hline

$y^2$  &   $p/(x+m)$,\quad  $ \text{$p=0, -15, \frac{-10}{3},
\frac{-5}{3}$
}$ & 2& \\[4pt]
     & As given by  $E_2=0$             & 2  &      \\[1pt]
     & As given by  \eqref{eq:c0y2dim1}, \; with $\theta=0$ &1&\\[1pt]\hline
$y^2+ \theta,\; \theta \neq 0$  &  $5 (\sqrt{\theta}i/3)^{1/2} \tan((\sqrt{\theta}i/3)^{1/2} (x+m))$ & 1 & \\[4pt]
                  & As given by  $E_2=0$               & 2  &      \\[1pt]
                 & As given by  \eqref{eq:c0y2dim1}  &1&\\[1pt]\hline
$y^{-1}$& $M, \; M\neq0$  & 2  &  \\[4pt]
        & $M/x,\quad  M\neq0$ & 1   & $x\pd_x + y \pd_y$    \\[4pt]\hline
 $y^{-3}$ & 0  & 3  &    \\
          & $M/x,\quad $ $M\neq 0$ & 1  & $2x \pd_x + y\pd_y$   \\
          & As given by  \eqref{eq:c5yndim1}, with $n=-3$  &1   &      \\[1pt]\hline
\parbox[c]{1.4in}{$ y^{\,n}$,\;\\ $(n\neq -3) $}& $\frac{-(n+3)/x}{(n+1)}, \; n\neq -1$ &2   & $(k_2 + k_1 x)\pd_x - \frac{2k_1}{n-1}\pd_y $   \\
                                                         &0  & 2  &    \\
     &$M/x,\quad  M\neq0, -\frac{n+3}{n+1},\; n\neq -1$  & 1  & $(n-1)x \pd_x -
2y \pd_y$   \\
& As given by  \eqref{eq:c5yndim1}                 &1   &      \\[1pt]\hline
$y^{-1} + \lambda y,\;\lambda \neq0$& $ \lambda x + m$  & 2  &    \\[4pt]
        & As given by  \eqref{eq:c5nzdim1}, with $n=-1$                 &1   &      \\[1pt]\hline
 $y^{-3}+ \lambda y,\;\lambda \neq0$ & 0  & 3  &    \\
          & As given by \eqref{eq:condam3}& 1   &    \\\hline
\parbox[c]{1.4in}{$ y^{\,n}+ \lambda y,\;\lambda \neq0$,\;\\ $(n\neq -3, -1) $}& $\frac{(3+n) \sqrt{\lambda } } { \sqrt{2(1+n)}} \tan
\left[\frac{\sqrt{\lambda(n+1)}\left(x+2
(3+n)m\right)}{\sqrt{2}}\right]$ &2   &    \\
& As given by  \eqref{eq:c5nzdim1}                 &1   &      \\[1pt]\hline

 $ F(y)$            & $M$&  1 & $  \pd_x$   \\\hline\hline
\end{tabular}
}
\end{table}
\subsection{Case 2: $F''=0$} By a linearization test for second-order
ODEs due to Lie \cite{lielint}, \eqref{eq:m} is linearizable for a
given function $F$ if and only if $F$ is linear,  and we are
therefore in the linear case of \eqref{eq:m}. It then follows from
an already cited result of Lie \cite{lie-ovsi1} that \eqref{eq:m}
has a symmetry algebra of maximal dimension eight. \par
We now undertake to verify that for linear $F= \lambda y + \theta,$
the dimension of $L$ is always eight, regardless of the value of
$A.$ It follows from the equivalence relations on functions $F$ that
we may assume that $F= \lambda y, \; \lambda \neq 0$ or $F= \theta.$
For $F= \lambda y,$ a substitution of this expression into
\eqref{eq:red1detm} shows that
\begin{equation}\label{eq:siglin}
\sigma = k_1+ \int \frac{1}{2} (A'\beta + A \beta' + \beta'') dx.
\end{equation}
When the latter expression for $\sigma$ is also substituted into
\eqref{eq:red1detm}, we obtain after a split into powers of $y$ the
following equations
\begin{subequations}\label{eq:detlin1}
\begin{align}
 \alpha  \left(-\lambda +A'\right)+A \alpha'+\alpha''&=0 \label{eq:detlin1a}\\
 - \lambda  \tau - A \tau\,'+ \tau\,'' &=0 \label{eq:detlin1b}\\
  A \alpha  \lambda - A \alpha  A'- A^2 \alpha'+ \left(\left(-\lambda + A'\right) \alpha'+\alpha  A''+\alpha\,'''\right)&=0 \label{eq:detlin1c}\\
 -A\, \beta\,  A'-A^2 \beta\,'+2 \left(-2 \lambda +A'\right) \beta\,'+\beta  A''+
 \beta\,'''&=0 \label{eq:detlin1d}
\end{align}
\end{subequations}
We note that only \eqref{eq:detlin1a} and \eqref{eq:detlin1c} are
dependent equations, and  we denote by $E_7$ and $E_8$ the left hand
side of \eqref{eq:detlin1c} and \eqref{eq:detlin1a}, respectively.
We have $E_7= E_8'- A E_8,$ where $E_8'= d E_8/dx,$ and this shows
that \eqref{eq:detlin1a} alone is equivalent to the system
consisting of the two equations \eqref{eq:detlin1a} and
\eqref{eq:detlin1c}. Therefore, to find $\alpha,$ $\tau$ and
$\beta,$ we only need to solve the independent linear equations
\eqref{eq:detlin1a}, \eqref{eq:detlin1b}, and \eqref{eq:detlin1d}.
The sum of  orders of these independent equations together with the
arbitrary constant $k_1$ in \eqref{eq:siglin} is exactly eight and
therefore yields  an $8$-parameter symmetry algebra, regardless of
the value of $A$. However, we cannot find explicit expressions for
the generator $V$ for arbitrary values of $A$ and $\lambda,$ because
this involves solving equations of the form \eqref{eq:detlin1b} for
which the general solution is not available. A similar analysis
holds when $F= \theta$ is a constant function.

\subsection{Summary}
 We have represented  the classification
results for the nonlinear case of \eqref{eq:m} in Table
\ref{tb:summary}, in which the first column indicates admissible
canonical forms of $F$ for which symmetries exist, while the second
column indicates the corresponding expression of the function $A.$
When  $A$ is given by complicated nonlinear or linear equations for
which the solution is not available, the required expression is
replaced in the column the determining equation for $A$ given in the
text. The third column gives the generator $V$ of $L$, but only for
those cases for which the explicit expression for $V$ is available,
and when its size is sufficiently small to fit in the table.
However,  explicit forms of the symmetry generator $V$ from which
symmetries can be readily calculated for given values of $A$ is
provided every time in the text for every possible canonical form
$F,$ and most of the generators $V$ with a relatively prominent size
 have been determined in the text whenever $A$ was known.\par

In the last row, the pair $\set{M, F(y)}$ represents an equation
with an arbitrary admissible function $F(y)$ and with $A=M,$
provided that such a function is not represented else where in the
table. For instance for $F(y)= y^{-1}$ and $A=M,$ $M\neq0,$ the
equation is represented in the row with $F=y^{-1},$ while for $A=0,$
the corresponding equation is represented in the block of rows with
$F= y^{n},\; (n\neq -3).$ Indeed, table rows represent non
equivalent equations.

\section{concluding remarks}
In this paper, we have given a group classification of equations of
the form \eqref{eq:m} and shown that the only admissible canonical
forms of $F$ admitting  symmetries are given by functions listed in
\eqref{eq:admis}. Moreover, we have shown that  in the nonlinear
case, the maximal dimension of the symmetry algebra $L$ achieved is
three, which coincide with the  result obtained in \cite{ovsi1} for
equations of the form \eqref{eq:ovsi1}.  We have also been able  to
determine for any canonical form $F$ explicit expressions of $A$ for
which $L$ has a given dimension $n$, where $0\leq n \leq 3, $ with a
few exceptions which apply mostly for cases of one dimensional
subalgebras, and rarely for two-dimensional subalgebras. Indeed, we
have not been able to provide explicit expressions of $A$ for which
$\dim L=2$ only for $F= \mu  e^y + \theta,$ and for $F= y^2 +
\theta$ although in the latter case  we have given some
one-parameter families of functions $A$ for which the equality $\dim
L=2$ occurs. These difficulties are due to the fact that  general
solutions of the related determining equations for $A,$ such as
$E_2=0,$ where $E_2$ is given by \eqref{eq:condac1c}, is not
available.

  For cases one one-dimensional symmetry algebras $L,$ the difficulty
with finding an explicit expression for the corresponding function
$A$ is often due to the fact that determining equations are quite
often complicated integro-differential equations or nonlinear
equations of the form \eqref{eq:condam3}. On the other hand, this
difficulty as well as the impossibility to give an explicit
expression of the symmetry generator in the linear case for
arbitrary values of $A$ is due to the fact that the general solution
is not known for  linear equations similar to \eqref{eq:detlin1b}
and of the form
\begin{equation*}
\beta'' + f(x) \beta' + g(x) \beta=0,
\end{equation*}
in which $f$ is an arbitrary function and $g$ is either another
arbitrary function or a nonzero constant. We note however that
equations of the latter form can always be reduced to the most
common canonical form
\begin{equation*}
w'' + h(x) w=0, \quad  \text{ with } h= -(f^2 - 4 g + 2f')/4,
\end{equation*}
by a change of the dependent variable of the form $\beta = w
e^{-(1/2) \int f dx}.$ Although the reduced equation depends on
fewer arbitrary functions, the difficulty with solving it remains
essentially the same for arbitrary functions $h.$ However, using the
provided determining equations, this classification can nevertheless
readily provide the explicit symmetry generator $V$ for any given
pair $\set{A, F}$ of labeling functions, or tell when a symmetry
does not exist.\par


%
\end{document}